\documentclass[10pt]{amsart}

\usepackage{amstext}
\usepackage{amsthm}
\usepackage{amssymb}

\makeatletter
\numberwithin{equation}{section}
\numberwithin{figure}{section}
\theoremstyle{plain}
\newtheorem{thm}{Theorem}
  \theoremstyle{plain}
  \newtheorem{prop}[thm]{Proposition}
  \theoremstyle{plain}
  \newtheorem{cor}[thm]{Corollary}
  \theoremstyle{remark}
  \newtheorem*{rem*}{Remark}

\allowdisplaybreaks
\usepackage[sort&compress,square,numbers]{natbib}

\makeatother

\begin{document}

\title{Convolved Fibonacci numbers and their applications}

\author{Taekyun Kim}

\address{Department of Mathematics, College of Science, Tianjin Polytechnic University,
Tianjin City, 300387, China\\ Department of Mathematics, Kwangwoon University, Seoul 139-701, Republic
of Korea}
\email{tkkim@kw.ac.kr}

\author{Dmitry V. Dolgy}

\address{Institute of Mathematics and Computer Science, Far Eastern Federal University, 690950 Vladivostok, Russia}

\email{dvdolgy@gmail.com}

\author{Dae San Kim}
\address{Department of Mathematics, Sogang University, Seoul 121-742, Republic
of Korea}

\email{dskim@sogang.ac.kr}

\author{Jong Jin Seo}

\address{Department of Applied Mathematics, Pukyong National University, Busan, Republic of Korea}

\email{seo2011@pknu.ac.kr}

\begin{abstract}
In this paper, we present a new approach to the convolved
Fibonacci numbers arising from the generating function of them and give
some new and explicit identities for the convolved Fibonacci numbers.
\end{abstract}

\keywords{Fibonacci numbers, convolved Fibonacci numbers, linear differential equation}
\subjclass[2010]{05A19, 11B83, 34A30}

\maketitle
\global\long\def\relphantom#1{\mathrel{\phantom{{#1}}}}

\section{Introduction}

As is well known, the Fibonacci numbers are given by
the numbers in the following integer sequnce:
\[
1,1,2,3,5,8,13,21,34,55,89,144,\dots.
\]

The sequence $F_{n}$ of Fibonacci numbers is defined by the recurrence
relation as follows:
\begin{equation}
F_{0}=1,\quad F_{1}=1,\quad F_{n}=F_{n-1}+F_{n-2},\quad\left(n\ge2\right),\quad\left(\text{see \cite{key-1,key-2,key-3,key-4,key-5,key-6,key-7,key-8}}\right).\label{eq:1}
\end{equation}

The sequence can be extended to negative index $n$ arising from the
re-arranged recurrence relation
\begin{equation}
F_{n-2}=F_{n}-F_{n-1},\quad\left(\text{see \cite{key-1,key-2,key-3,key-4,key-5,key-6,key-7,key-8,key-9,key-10,key-11,key-12,key-13}}\right)\label{eq:2}
\end{equation}
which yields the sequence of ``negafibonacci'' numbers satisfying
\begin{equation}
F_{-n}=\left(-1\right)^{n+1}F_{n},\quad\left(\text{see \cite{key-11,key-12}}\right).\label{eq:3}
\end{equation}

It is well known that the generating function of Fibonacci numbers
is given by
\begin{equation}
\frac{1}{1-t-t^{2}}=\sum_{n=0}^{\infty}F_{n}t^{n},\quad\left(\text{see \cite{key-3,key-4,key-5,key-6}}\right).\label{eq:4}
\end{equation}

The convolved Fibonacci numbers $p_{n}\left(x\right)$, $\left(n\ge0\right)$,
are defined by the generating function
\begin{equation}
\left(\frac{1}{1-t-t^{2}}\right)^{x}=\sum_{n=0}^{\infty}p_{n}\left(x\right)\frac{t^{n}}{n!},\quad\left(x\in\mathbb{R}\right).\label{eq:5}
\end{equation}

From (\ref{eq:4}) and (\ref{eq:5}), we note that
\begin{equation}
\frac{p_{n}\left(1\right)}{n!}=F_{n},\quad\left(n\ge0\right).\label{eq:6}
\end{equation}

In this paper, we present a new approach to the convolved
Fibonacci numbers arising from the generating function of them and give
some new and explicit identities for the convolved Fibonacci numbers.

\section{Convolved Fibonacci numbers and their applications}

From (\ref{eq:51}), we note that
\begin{align}
\sum_{n=0}^{\infty}p_{n}\left(x\right)\frac{t^{n}}{n!} & =\left(\frac{1}{1-t-t^{2}}\right)^{x}=\left(\frac{1}{1-t-t^{2}}\right)\left(\frac{1}{1-t-t^{2}}\right)^{x-1}\label{eq:7}\\
 & =\left(\sum_{l=0}^{\infty}p_{l}\left(1\right)\frac{t^{l}}{l!}\right)\left(\sum_{m=0}^{\infty}p_{m}\left(x-1\right)\frac{t^{m}}{m!}\right)\nonumber \\
 & =\sum_{n=0}^{\infty}\left(\sum_{l=0}^{n}\binom{n}{l}p_{l}\left(1\right)p_{n-l}\left(x-1\right)\right)\frac{t^{n}}{n!}.\nonumber
\end{align}

By comparing the coefficients on both sides of (\ref{eq:7}),
we obtain the following proposition.
\begin{prop}
\label{prop:1}For $n\ge0$, $x\in\mathbb{R}$, we have
\[
p_{n}\left(x\right)=\sum_{l=0}^{n}\binom{n}{l}p_{l}\left(1\right)p_{n-l}\left(x-1\right).
\]

\end{prop}
Let us take $x=r\in\mathbb{N}$. Then, by Proposition \ref{prop:1},
we get
\begin{align}
p_{n}\left(r\right) & =\sum_{l_{1}=0}^{n}\binom{n}{l_{1}}p_{l_{1}}\left(1\right)p_{n-l_{1}}\left(r-1\right)\label{eq:8}\\
 & =\sum_{l_{1}=0}^{n}\sum_{l_{2}=0}^{n-l_{1}}\binom{n}{l_{1}}\binom{n-l_{1}}{l_{2}}p_{l_{1}}\left(1\right)p_{l_{2}}\left(1\right)p_{n-l_{1}-l_{2}}\left(r-2\right)\nonumber \\
 & =\sum_{l_{1}=0}^{n}\sum_{l_{2}=0}^{n-l_{1}}\sum_{l_{3}=0}^{n-l_{1}-l_{2}}\binom{n}{l_{1}}\binom{n-l_{1}}{l_{2}}\binom{n-l_{1}-l_{2}}{l_{3}}p_{l_{1}}\left(1\right)p_{l_{2}}\left(1\right)\nonumber\\
 &\relphantom{=}\times p_{l_{3}}\left(1\right)p_{n-l_{1}-l_{2}-l_{3}}\left(r-3\right)\nonumber \\
 & \vdots\nonumber \\
 & =\sum_{l_{1}=0}^{n}\sum_{l_{2}=0}^{n-l_{1}}\cdots\sum_{l_{r-1}=0}^{n-l_{1}-\cdots-l_{r-2}}\binom{n}{l_{1}}\binom{n-l_{1}}{l_{2}}\cdots\nonumber\\
 &\relphantom{=}\times \binom{n-l_{1}-l_{2}-\cdots-l_{r-2}}{l_{r-1}}\nonumber\\
 &\times \left(\prod_{k=1}^{r-1}p_{l_{k}}\left(1\right)\right)p_{n-l_{1}-l_{2}-\cdots-l_{r-1}}\left(1\right).\nonumber
\end{align}

Therefore, by (\ref{eq:8}), we obtain the following corollary.
\begin{cor}
\label{cor:2} For $r\in\mathbb{N}$ and $n\ge0$, we have
\begin{align*}
p_{n}\left(r\right)&=\sum_{l_{1}=0}^{n}\sum_{l_{2}=0}^{n-l_{1}}\cdots\sum_{l_{r-1}=0}^{n-l_{1}-\cdots-l_{r-2}}\binom{n}{l_{1}}\binom{n-l_{1}}{l_{2}}\cdots\\
&\relphantom{=}\times\binom{n-l_{1}-l_{2}-\cdots-l_{r-2}}{l_{r-1}}\left(\prod_{k=1}^{r-1}p_{l_{k}}\left(1\right)\right)p_{n-l_{1}-l_{2}-\cdots-l_{r-1}}\left(1\right).
\end{align*}

\end{cor}
We observe that
\begin{align}
\left(\frac{1}{1-t-t^{2}}\right)^{x} & =\left(\frac{1}{1-t-t^{2}}\right)^{r}\left(\frac{1}{1-t-t^{2}}\right)^{x-r}\label{eq:9}\\
 & =\left(\sum_{l=0}^{\infty}p_{l}\left(r\right)\frac{t^{l}}{l!}\right)\left(\sum_{m=0}^{\infty}p_{m}\left(x-r\right)\frac{t^{m}}{m!}\right)\nonumber \\
 & =\sum_{n=0}^{\infty}\left(\sum_{l=0}^{n}\binom{n}{l}p_{l}\left(r\right)p_{n-l}\left(x-r\right)\right)\frac{t^{n}}{n!}.\nonumber
\end{align}

Therefore, by (\ref{eq:51}) and (\ref{eq:9}), we obtain the following
theorem.
\begin{thm}
\label{thm:3} For $n\ge0$, $r\in\mathbb{N}$, we have
\[
p_{n}\left(x\right)=\sum_{l=0}^{n}\binom{n}{l}p_{l}\left(r\right)p_{n-l}\left(x-r\right)=\sum_{l=0}^{n}\binom{n}{l}p_{n-l}\left(r\right)p_{l}\left(x-r\right).
\]

\end{thm}
Let us take $x=r+1$ in Theorem \ref{thm:3}. Then, we have
\begin{align}
p_{n}\left(r+1\right) & =\sum_{l=0}^{n}\binom{n}{l}p_{n-l}\left(r\right)p_{l}\left(1\right)\label{eq:10}\\
 & =\sum_{l=0}^{n}\binom{n}{l}p_{n-l}\left(r\right)l!\frac{p_{l}\left(1\right)}{l!}\nonumber \\
 & =\sum_{l=0}^{n}\left(n\right)_{l}p_{n-l}\left(r\right)F_{l},\nonumber
\end{align}
where $\left(x\right)_{n}=x\left(x-1\right)\cdots\left(x-n+1\right)$,
$\left(n\ge1\right)$, $\left(x\right)_{0}=1$.
\begin{cor}
\label{cor:4} For $r\in\mathbb{N}$, $n\ge0$, we have
\[
p_{n}\left(r+1\right)=\sum_{l=0}^{n}\left(n\right)_{l}p_{n-l}\left(r\right)F_{l}.
\]

\end{cor}
Taking $r=1$ in Corollary \ref{cor:4}, we have
\begin{align}
p_{n}\left(2\right) & =\sum_{l=0}^{n}\left(n\right)_{l}p_{n-l}\left(1\right)F_{l}\label{eq:11}\\
 & =\sum_{l=0}^{n}\left(n\right)_{l}\left(n-l\right)!\frac{p_{n-l}\left(1\right)}{\left(n-l\right)!}F_{l}\nonumber \\
 & =n!\sum_{n=0}^{n}\binom{n}{l}\binom{n}{l}^{-1}F_{n-l}F_{l}\nonumber \\
 & =n!\sum_{l=0}^{n}F_{l}F_{n-l}.\nonumber
\end{align}

Thus, by (\ref{eq:11}), we get
\begin{equation}
\frac{p_{n}\left(2\right)}{n!}=\sum_{l=0}^{n}F_{l}F_{n-l},\quad\left(n\ge0\right).\label{eq:12}
\end{equation}

From (\ref{eq:12}), we can also derive the following equation.

\begin{align}
p_{n}\left(3\right) & =\sum_{l_{1}=0}^{n}\left(n\right)_{l_{1}}p_{n-l_{1}}\left(2\right)F_{l_{1}}\label{eq:13}\\
 & =\sum_{l_{1}=0}^{n}\left(n\right)_{l_{1}}\left(n-l_{1}\right)!\frac{p_{n-l_{1}}\left(2\right)}{\left(n-l_{1}\right)!}F_{l_{1}}\nonumber \\
 & =n!\sum_{l_{1}=0}^{n}\sum_{l_{2}=0}^{n-l_{1}}F_{l_{1}}F_{l_{2}}F_{n-l_{1}-l_{2}}.\nonumber
\end{align}

Thus, by (\ref{eq:13}), we get
\begin{equation}
\frac{p_{n}\left(3\right)}{n!}=\sum_{l_{1}=0}^{n}\sum_{l_{2}=0}^{n-l_{1}}F_{l_{1}}F_{l_{2}}F_{n-l_{1}-l_{2}}.\label{eq:14}
\end{equation}

For $r=3$ in Corollary \ref{cor:4}, we have
\begin{align}
p_{n}\left(4\right) & =\sum_{l_{1}=0}^{n}\left(n\right)_{l_{1}}p_{n-l_{1}}\left(3\right)F_{l_{1}}\label{eq:15}\\
 & =n!\sum_{l_{1}=0}^{n}\frac{p_{n-l_{1}}\left(3\right)}{\left(n-l_{1}\right)!}F_{l_{1}}\nonumber \\
 & =n!\sum_{l_{1}=0}^{n}\sum_{l_{2}=0}^{n-l_{1}}\sum_{l_{3}=0}^{n-l_{1}-l_{2}}F_{l_{1}}F_{l_{2}}F_{l_{3}}F_{n-l_{1}-l_{2}-l_{3}}.\nonumber
\end{align}

From (\ref{eq:15}), we note that
\begin{equation}
\frac{p_{n}\left(4\right)}{n!}=\sum_{l_{1}=0}^{n}\sum_{l_{2}=0}^{n-l_{1}}\sum_{l_{3}=0}^{n-l_{1}-l_{2}}F_{l_{1}}F_{l_{2}}F_{l_{3}}F_{n-l_{1}-l_{2}-l_{3}}.\label{eq:16}
\end{equation}

Continuing this process, we have
\begin{equation}
\frac{p_{n}\left(r+1\right)}{n!}=\sum_{l_{1}=0}^{n}\sum_{l_{2}=0}^{n-l_{1}}\cdots\sum_{l_{r}=0}^{n-l_{1}-\cdots-l_{r-1}}F_{l_{1}}F_{l_{2}}\cdots F_{l_{r}}F_{n-l_{1}-l_{2}-\cdots-l_{r}},\label{eq:17}
\end{equation}
where $r\in\mathbb{N}$.
\begin{thm}
\label{thm:5} For $r\in\mathbb{N}$ and $n\ge0$, we have
\[
\frac{p_{n}\left(r+1\right)}{n!}=\sum_{l_{1}=0}^{n}\sum_{l_{2}=0}^{n-l_{1}}\cdots\sum_{l_{r}=0}^{n-l_{1}-\cdots-l_{r-1}}F_{l_{1}}F_{l_{2}}\cdots F_{l_{r}}F_{n-l_{1}-l_{2}-\cdots-l_{r}}.
\]

\end{thm}
Let
\begin{align}
F\left(t,x\right) & =\left(1-t-t^{2}\right)^{-x}\label{eq:18}\\
 & =e^{-x\log\left(1-t-t^{2}\right)}.\nonumber
\end{align}

Then, by (\ref{eq:18}), we get
\begin{align}
F^{\left(1\right)}\left(t,x\right) & =\frac{dF}{dt}\left(t,x\right)\label{eq:19}\\
 & =x\left(1+2t\right)\left(1-t-t^{2}\right)^{-x-1}\nonumber \\
 & =x\left(1+2t\right)F\left(t,x+1\right),\nonumber
\end{align}
\begin{align}
F^{\left(2\right)}\left(t,x\right) & =\frac{dF^{\left(1\right)}}{dt}\left(t,x\right)\label{eq:20}\\
 & =2xF\left(t,x+1\right)+\left\langle x\right\rangle _{2}\left(1+2t\right)^{2}F\left(t,x+2\right),\nonumber
\end{align}
where $\left\langle x\right\rangle _{n}=x\left(x+1\right)\cdots\left(x+n-1\right)$,
$\left(n\ge1\right)$, $\left\langle x\right\rangle _{0}=1$.

From (\ref{eq:20}), we note that
\begin{align}
F^{\left(3\right)}\left(t,x\right) & =\frac{dF^{\left(2\right)}}{dt}\left(t,x\right)\label{eq:21}\\
 & =6\left\langle x\right\rangle _{2}\left(1+2t\right)F\left(t,x+2\right)+\left\langle x\right\rangle _{3}\left(1+2t\right)^{3}F\left(t,x+3\right).\nonumber
\end{align}
\begin{align}
F^{\left(4\right)}\left(t,x\right) & =\frac{dF^{\left(3\right)}}{dt}\left(t,x\right)\label{eq:22}\\
 & =12\left\langle x\right\rangle _{2}F\left(t,x+2\right)+12\left\langle x\right\rangle _{3}\left(1+2t\right)^{2}F\left(t,x+3\right)\nonumber \\
 & \relphantom =+\left\langle x\right\rangle _{4}\left(1+2t\right)^{4}F\left(t,x+4\right),\nonumber
\end{align}
\begin{align}
F^{\left(5\right)}\left(t,x\right) & =\frac{dF^{\left(4\right)}}{dt}\left(t,x\right)\label{eq:23}\\
 & =60\left\langle x\right\rangle _{3}\left(1+2t\right)F\left(t,x+3\right)+20\left\langle x\right\rangle _{4}\left(1+2t\right)^{3}F\left(t,x+4\right)\nonumber\\
 &\relphantom{=}+\left\langle x\right\rangle _{5}\left(1+2t\right)^{5}F\left(t,x+5\right)\nonumber
\end{align}
 and
\begin{align}
F^{\left(6\right)}\left(t,x\right) & =\frac{dF^{\left(5\right)}}{dt}\left(t,x\right)\label{eq:24}\\
 & =120\left\langle x\right\rangle _{3}F\left(t,x+3\right)+180\left\langle x\right\rangle _{4}\left(1+2t\right)^{2}F\left(t,x+4\right)\nonumber \\
 & \relphantom =+30\left\langle x\right\rangle _{5}\left(1+2t\right)^{4}F\left(t,x+5\right)+\left\langle x\right\rangle _{6}\left(1+2t\right)^{6}F\left(t,x+6\right).\nonumber
\end{align}

Thus, we are led to put
\begin{align}
&\relphantom{=}F^{\left(N\right)}\left(t,x\right)=\left(\frac{d}{dt}\right)^{N}F\left(t,x\right)\label{eq:25}\\
&=\sum_{i=0}^{\left[\frac{N}{2}\right]}a_{i}\left(N\right)\left\langle x\right\rangle _{N-i}\left(1+2t\right)^{N-2i}F\left(t,x+N-i\right)\nonumber
\end{align}
where $N\in\mathbb{N}$.

Taking the derivatives of (\ref{eq:25}) with respect to $t$, we have
\begin{align}
F^{\left(N+1\right)}\left(t,x\right) & =\sum_{i=0}^{\left[\frac{N}{2}\right]}a_{i}\left(N\right)\left\langle x\right\rangle _{N-i}\left(1+2t\right)^{N-2i}F^{\left(1\right)}\left(t,x+N-i\right)\label{eq:26}\\
 & \relphantom =+\sum_{i=0}^{\left[\frac{N}{2}\right]}a_{i}\left(N\right)\left\langle x\right\rangle _{N-i}2\left(N-2i\right)\left(1+2t\right)^{N-2i-1}F\left(t,x+N-i\right)\nonumber \\
 & =\sum_{i=0}^{\left[\frac{N}{2}\right]}2\left(N-2i\right)a_{i}\left(N\right)\left\langle x\right\rangle _{N-i}\left(1+2t\right)^{N-2i-1}F\left(t,x+N-i\right)\nonumber \\
 & \relphantom =+\sum_{i=0}^{\left[\frac{N}{2}\right]}a_{i}\left(N\right)\left\langle x\right\rangle _{N-i+1}\left(1+2t\right)^{N-2i+1}F\left(t,x+N-i+1\right)\nonumber \\
 & =\sum_{i=1}^{\left[\frac{N}{2}\right]+1}2\left(N-2i+2\right)a_{i-1}\left(N\right)\left\langle x\right\rangle _{N-i+1}\nonumber\\
 &\relphantom{=}\times\left(1+2t\right)^{N-2i+1}F\left(t,x+N-i+1\right)\nonumber \\
 & \relphantom =+\sum_{i=0}^{\left[\frac{N}{2}\right]}a_{i}\left(N\right)\left\langle x\right\rangle _{N-i+1}\left(1+2t\right)^{N-2i+1}F\left(t,x+N-i+1\right).\nonumber
\end{align}

On the other hand, by replacing $N$ by $N+1$ in (\ref{eq:25}),
we get
\begin{equation}
F^{\left(N+1\right)}\left(t,x\right)=\sum_{i=0}^{\left[\frac{N+1}{2}\right]}a_{i}\left(N+1\right)\left\langle x\right\rangle _{N-i+1}\left(1+2t\right)^{N-2i+1}F\left(t,x+N-i+1\right).\label{eq:27}
\end{equation}

\emph{Case 1. }Let $N$ be an even number. Then we have
\begin{align}
 & \sum_{i=1}^{\frac{N}{2}+1}2\left(N-2i+2\right)a_{i-1}\left(N\right)\left\langle x\right\rangle _{N-i+1}\left(1+2t\right)^{N-2i+1}F\left(t,x+N-i+1\right)\label{eq:28}\\
 & +\sum_{i=0}^{\frac{N}{2}}a_{i}\left(N\right)\left\langle x\right\rangle _{N-i+1}\left(1+2t\right)^{N-2i+1}F\left(t,x+N-i+1\right)\nonumber \\
 & =\sum_{i=0}^{\frac{N}{2}}a_{i}\left(N+1\right)\left\langle x\right\rangle _{N-i+1}\left(1+2t\right)^{N-2i+1}F\left(t,x+N-i+1\right).\nonumber
\end{align}

Comparing the coefficients on both sides of (\ref{eq:28}), we
get
\begin{align}
a_{0}\left(N+1\right) & =a_{0}\left(N\right),\label{eq:29}\\
a_{i}\left(N+1\right) & =2\left(N-2i+2\right)a_{i-1}\left(N\right)+a_{i}\left(N\right),\quad\left(1\le i\le\frac{N}{2}\right).\label{eq:30}
\end{align}

\emph{Case 2}. Let $N$ be an odd number. Then, we have
\begin{align}
 & \sum_{i=1}^{\frac{N+1}{2}}2\left(N-2i+2\right)a_{i-1}\left(N\right)\left\langle x\right\rangle _{N-i+1}\left(1+2t\right)^{N-2i+1}F\left(t,x+N-i+1\right)\label{eq:31}\\
 & +\sum_{i=0}^{\frac{N-1}{2}}a_{i}\left(N\right)\left\langle x\right\rangle _{N-i+1}\left(1+2t\right)^{N-2i+1}F\left(t,x+N-i+1\right)\nonumber \\
 & =\sum_{i=0}^{\frac{N+1}{2}}a_{i}\left(N+1\right)\left\langle x\right\rangle _{N-i+1}\left(1+2t\right)^{N-2i+1}F\left(t,x+N-i+1\right).\nonumber
\end{align}

Comparing the coefficients on both sides of (\ref{eq:31}), we
have
\begin{equation}
a_{0}\left(N+1\right)=a_{0}\left(N\right),\quad a_{\frac{N+1}{2}}\left(N+1\right)=2a_{\frac{N-1}{2}}\left(N\right),\label{eq:32}
\end{equation}
and
\begin{equation}
a_{i}\left(N+1\right)=2\left(N-2i+2\right)a_{i-1}\left(N\right)+a_{i}\left(N\right),\quad\left(1\le i\le\frac{N-1}{2}\right).\label{eq:33}
\end{equation}

In addition, we have the following ``initial conditions'':
\begin{equation}
F^{\left(0\right)}\left(t,x\right)=F\left(t,x\right)=a_{0}\left(0\right)F\left(t,x\right).\label{eq:34}
\end{equation}

Thus, by (\ref{eq:34}), we get $a_{0}\left(0\right)=1$.

From (\ref{eq:19}) and (\ref{eq:25}), we note that
\begin{equation}
F^{\left(1\right)}\left(t,x\right)=a_{0}\left(1\right)x\left(1+2t\right)F\left(t,x+1\right)=x\left(1+2t\right)F\left(t,x+1\right).\label{eq:35}
\end{equation}

Thus, by (\ref{eq:35}), we see that $a_{0}\left(1\right)=1$.

By (\ref{eq:20}) and (\ref{eq:25}), we easily get

\begin{align}
F^{\left(2\right)}\left(t,x\right) & =\sum_{i=0}^{1}a_{i}\left(2\right)\left\langle x\right\rangle _{2-i}\left(1+2t\right)^{2-2i}F\left(t,x+2-i\right)\label{eq:36}\\
 & =a_{0}\left(2\right)\left\langle x\right\rangle _{2}\left(1+2t\right)^{2}F\left(t,x+2\right)+a_{1}\left(2\right)xF\left(t,x+1\right)\nonumber \\
 & =2xF\left(t,x+1\right)+\left\langle x\right\rangle _{2}\left(1+2t\right)^{2}F\left(t,x+2\right).\nonumber
\end{align}

Thus, by comparing the coefficients on both sides of (\ref{eq:36}),
we get
\begin{equation}
a_{0}\left(2\right)=1,\quad\text{and}\quad a_{1}\left(2\right)=2.\label{eq:37}
\end{equation}

In (\ref{eq:25}), it is not difficult to show that
\begin{equation}
a_{\frac{N+1}{2}}\left(N\right)=0,\quad\text{for all positive integers }N.\label{eq:38}
\end{equation}

From (\ref{eq:38}), we note that
\begin{equation}
a_{1}\left(1\right)=a_{2}\left(3\right)=a_{3}\left(5\right)=a_{4}\left(7\right)=\cdots=0.\label{eq:39}
\end{equation}

By (\ref{eq:38}), we get
\begin{equation}
F^{\left(N\right)}\left(t,x\right)=\sum_{i=0}^{\left[\frac{N+1}{2}\right]}a_{i}\left(N\right)\left\langle x\right\rangle _{N-i}\left(1+2t\right)^{N-2i}F\left(t,x+N-i\right),\label{eq:40}
\end{equation}
where
\begin{equation}
a_{0}\left(N+1\right)=a_{0}\left(N\right),\quad a_{\frac{N+1}{2}}\left(N\right)=0,\quad\text{for all positive integers }N,\label{eq:41}
\end{equation}
and
\begin{equation}
a_{i}\left(N+1\right)=2\left(N-2i+2\right)a_{i-1}\left(N\right)+a_{i}\left(N\right),\quad\left(1\le i\le\left[\frac{N+1}{2}\right]\right).\label{eq:42}
\end{equation}

From (\ref{eq:41}), we note that
\begin{equation}
a_{0}\left(N+1\right)=a_{0}\left(N\right)=a_{0}\left(N-1\right)=\cdots=a_{0}\left(1\right)=1.\label{eq:43}
\end{equation}

For $i=1,2,3$ in (\ref{eq:42}), we have
\begin{align}
a_{1}\left(N+1\right) & =2\sum_{k=0}^{N-1}\left(N-k\right)a_{0}\left(N-k\right),\label{eq:44}\\
a_{2}\left(N+1\right) & =2\sum_{k=0}^{N-3}\left(N-2-k\right)a_{1}\left(N-k\right),\label{eq:45}
\end{align}
and
\begin{equation}
a_{3}\left(N+1\right)=2\sum_{k=0}^{N-5}\left(N-4-k\right)a_{2}\left(N-k\right).\label{eq:46}
\end{equation}

Thus, we can deduce that, for $1\le i\le\left[\frac{N+1}{2}\right],$
\begin{align}
a_{i}\left(N+1\right) & =2\sum_{k=0}^{N-2i+1}\left(N-k-2i+2\right)a_{i-1}\left(N-k\right)\label{eq:47}\\
 & =2\sum_{k=0}^{N+2-2i}ka_{i-1}\left(k+2i-2\right).\nonumber
\end{align}

Now, we give explicit expressions for $a_{i}\left(N+1\right)$.

From (\ref{eq:43}), (\ref{eq:44}), (\ref{eq:45}) and (\ref{eq:46}),
we have
\begin{align}
a_{1}\left(N+1\right) & =2\sum_{k_{1}=1}^{N}k_{1}a_{0}\left(k_{1}\right)=2\sum_{k_{1}=1}^{N}k_{1}\label{eq:48}\\
a_{2}\left(N+1\right) & =2\sum_{k_{2}=1}^{N-2}k_{2}a_{1}\left(k_{2}+2\right)=2^{2}\sum_{k_{2}=1}^{N-2}\sum_{k_{1}=1}^{k_{2}+1}k_{2}k_{1},\label{eq:49}\\
a_{3}\left(N+1\right) & =2\sum_{k_{3}=1}^{N-4}k_{3}a_{2}\left(k_{3}+4\right)=2^{3}\sum_{k_{3}=1}^{N-4}\sum_{k_{2}=1}^{k_{3}+1}\sum_{k_{1}=1}^{k_{2}+1}k_{3}k_{2}k_{1}\label{eq:50}
\end{align}
and
\begin{equation}
a_{4}\left(N+1\right)=2^{4}\sum_{k_{4}=1}^{N-6}\sum_{k_{3}=1}^{k_{4}+1}\sum_{k_{2}=1}^{k_{3}+1}\sum_{k_{1}=1}^{k_{2}+1}k_{4}k_{3}k_{2}k_{1}.\label{eq:51}
\end{equation}

Continuing this process, we have
\begin{equation}
a_{i}\left(N+1\right)=2^{i}\sum_{k_{i}=1}^{N-2i+2}\sum_{k_{i-1}=1}^{k_{i}+1}\cdots\sum_{k_{1}=1}^{k_{2}+1}\left(\prod_{l=1}^{i}k_{l}\right),\quad\left(1\le i\le\left[\frac{N+1}{2}\right]\right).\label{eq:52}
\end{equation}

Therefore, by (\ref{eq:52}), we obtain the following theorem.
\begin{thm}
\label{thm:6} For $N=0,1,2,\dots$, the family of differential equations
\begin{align*}
F^{\left(N\right)}\left(t,x\right)&=\left(\frac{d}{dt}\right)^{N}F\left(t,x\right)\\
&=\left(\sum_{i=0}^{\left[\frac{N+1}{2}\right]}a_{i}\left(N\right)\left\langle x\right\rangle _{N-i}\left(1+2t\right)^{N-2i}\left(1-t-t^{2}\right)^{-N+i}\right)F\left(t,x\right)
\end{align*}
have a solution
\[
F\left(t,x\right)=\left(1-t-t^{2}\right)^{-x}
\]
where
\[
a_{0}\left(N\right)=1,\quad a_{\frac{N+1}{2}}\left(N\right)=0,\quad\text{for all positive integers }N,
\]
and
\[
a_{i}\left(N\right)=2^{i}\sum_{k_{i}=1}^{N-2i+1}\sum_{k_{i-1}=1}^{k_{i}+1}\cdots\sum_{k_{1}=1}^{k_{2}+1}\left(\prod_{l=1}^{i}k_{l}\right),\quad\left(1\le i\le\left[\frac{N}{2}\right]\right).
\]

\end{thm}
From (\ref{eq:4}), we note that
\begin{equation}
1=\sum_{k=0}^{\infty}F_{k}t^{k}\left(1-t-t^{2}\right)=\sum_{k=0}^{\infty}F_{k}t^{k}-\sum_{k=1}^{\infty}F_{k-1}t^{k}-\sum_{k=2}^{\infty}F_{k-2}t^{k}.\label{eq:53}
\end{equation}

Comparing the coefficients on the both sides of (\ref{eq:53}), we
have
\begin{equation}
F_{0}=1,\quad F_{1}-F_{0}=0\iff F_{1}=F_{0}=1,\label{eq:54}
\end{equation}
and
\begin{equation}
F_{k}-F_{k-1}-F_{k-2}=0\quad\text{if }k\ge2.\label{eq:55}
\end{equation}

By (\ref{eq:4}), we easily get
\begin{equation}
F\left(t,x\right)=\left(1-t-t^{2}\right)^{-x}=\sum_{k=0}^{\infty}p_{k}\left(x\right)\frac{t^{k}}{k!},\label{eq:56}
\end{equation}
and
\begin{equation}
F^{\left(N\right)}\left(t,x\right)=\left(\frac{d}{dt}\right)^{N}F\left(t,x\right)=\sum_{k=0}^{\infty}p_{k+N}\left(x\right)\frac{t^{k}}{k!}.\label{eq:57}
\end{equation}

On the other hand, by Theorem \ref{thm:6}, we get
\begin{align}
F^{\left(N\right)}\left(t,x\right) & =\sum_{i=0}^{\left[\frac{N+1}{2}\right]}a_{i}\left(N\right)\left\langle x\right\rangle _{N-i}\left(1+2t\right)^{N-2i}F\left(t,x+N-i\right)\label{eq:58}\\
 & =\sum_{i=0}^{\left[\frac{N+1}{2}\right]}a_{i}\left(N\right)\left\langle x\right\rangle _{N-i}\left(1+2t\right)^{N-2i}\sum_{m=0}^{\infty}p_{m}\left(x+N-i\right)\frac{t^{m}}{m!}\nonumber \\
 & =\sum_{i=0}^{\left[\frac{N+1}{2}\right]}a_{i}\left(N\right)\left\langle x\right\rangle _{N-i}\sum_{l=0}^{\infty}\left(N-2i\right)_{l}2^{l}\frac{t^{l}}{l!}\sum_{m=0}^{\infty}p_{m}\left(x+N-i\right)\frac{t^{m}}{m!}\nonumber \\
 & =\sum_{i=0}^{\left[\frac{N+1}{2}\right]}a_{i}\left(N\right)\left\langle x\right\rangle _{N-i}\sum_{k=0}^{\infty}\left(\sum_{l=0}^{k}\binom{k}{l}\left(N-2i\right)_{l}2^{l}p_{k-l}\left(x+N-i\right)\right)\frac{t^{k}}{k!}\nonumber \\
 & =\sum_{k=0}^{\infty}\left(\sum_{i=0}^{\left[\frac{N+1}{2}\right]}\sum_{l=0}^{k}\binom{k}{l}\left(N-2i\right)_{l}2^{l}a_{i}\left(N\right)\left\langle x\right\rangle _{N-i}p_{k-l}\left(x+N-i\right)\right)\frac{t^{k}}{k!}.\nonumber
\end{align}

Therefore, by comparing the coefficients on both sides of (\ref{eq:57})
and (\ref{eq:58}), we obtain the following theorem.
\begin{thm}
\label{thm:7} For $k,N=0,1,2,\dots$, we have
\[
p_{k+N}\left(x\right)=\sum_{i=0}^{\left[\frac{N+1}{2}\right]}\sum_{l=0}^{k}\binom{k}{l}\left(N-2i\right)_{l}2^{l}a_{i}\left(N\right)\left\langle x\right\rangle _{N-i}p_{k-l}\left(x+N-i\right),
\]
 where
\begin{align*}
a_{0}\left(N\right) & =1,\quad a_{\frac{N+1}{2}}\left(N\right)=0,\quad\text{for all positive integers }N,\\
a_{i}\left(N\right) & =2^{i}\sum_{k_{i}=1}^{N-2i+1}\sum_{k_{i-1}=1}^{k_{i}+1}\cdots\sum_{k_{1}=1}^{k_{2}+1}\left(\prod_{l=1}^{i}k_{l}\right),\quad\left(1\le i\le\left[\frac{N}{2}\right]\right).
\end{align*}

\end{thm}
When $k=0$ in Theorem \ref{thm:7}, we have the following corollary.
\begin{cor}
\label{cor:8} For $N=0,1,2,\dots$, we have
\[
p_{N}\left(x\right)=\sum_{i=0}^{\left[\frac{N+1}{2}\right]}a_{i}\left(N\right)\left\langle x\right\rangle _{N-i}.
\]

\end{cor}
Let us take $x=1$ in Corollary \ref{cor:8}. Then, we easily get
\begin{equation}
p_{N}\left(1\right)=\sum_{i=0}^{\left[\frac{N+1}{2}\right]}a_{i}\left(N\right)\left(N-i\right)!=N!+\sum_{i=1}^{\left[\frac{N+1}{2}\right]}a_{i}\left(N\right)\left(N-i\right)!\label{eq:59}
\end{equation}

Thus, by (\ref{eq:59}), we get
\begin{align}
\frac{p_{N}\left(1\right)}{N!} & =1+\frac{1}{N!}\sum_{i=1}^{\left[\frac{N+1}{2}\right]}a_{i}\left(N\right)\left(N-i\right)!\label{eq:60}\\
 & =1+\frac{1}{N!}\sum_{i=1}^{\left[\frac{N+1}{2}\right]}\sum_{k_{i}=1}^{N+1-2i}\sum_{k_{i-1}=1}^{k_{i}+1}\cdots\sum_{k_{1}=1}^{k_{2}+1}2^{i}\left(\prod_{l=1}^{i}k_{l}\right)\left(N-i\right)!\nonumber
\end{align}

Therefore, by (\ref{eq:6}) and (\ref{eq:60}), we obtain the following
corollary.
\begin{cor}
\label{cor:9}For $N=0,1,2,\dots$, we have
\[
F_{N}-1=\frac{1}{N!}\left(\sum_{i=1}^{\left[\frac{N+1}{2}\right]}\sum_{k_{i}=1}^{N+1-2i}\sum_{k_{i-1}=1}^{k_{i}+1}\cdots\sum_{k_{1}=1}^{k_{2}+1}2^{i}\left(\prod_{l=1}^{i}k_{l}\right)\left(N-i\right)!\right).
\]
\end{cor}
\begin{rem*}
Recently, several authors have studied special polynomials and sequences
arising from the generating functions (see \cite{key-1,key-2,key-3,key-4,key-5,key-6,key-7,key-8,key-9,key-10,key-11,key-12,key-13,key-14,key-15,key-16}).
\end{rem*}

\subsection*{Acknowledgements}
This paper is supported by grant NO 14-11-00022 of Russian Scientific Fund.

\bibliographystyle{amsplain}
\providecommand{\bysame}{\leavevmode\hbox to3em{\hrulefill}\thinspace}
\providecommand{\MR}{\relax\ifhmode\unskip\space\fi MR }
\providecommand{\MRhref}[2]{%
  \href{http://www.ams.org/mathscinet-getitem?mr=#1}{#2}
}
\providecommand{\href}[2]{#2}

\end{document}